\documentclass[11pt]{article}
\usepackage{hyperref}
\usepackage{amsmath,amssymb,amsthm}
\usepackage{bbm}
\usepackage{latexsym}
\usepackage{fullpage}
\usepackage{graphics}

\theoremstyle{definition}
\newtheorem{theorem}{Theorem}[section]

\DeclareMathOperator{\Exp}{E}
\DeclareMathOperator{\Var}{Var}

\newcommand{\set}[1]{\{#1\}}

\title{Stirling's Approximation for Central Extended Binomial Coefficients}

\author{Steffen Eger
}

\date{}

\begin{document}
\maketitle
\begin{abstract}
  We derive asymptotic formulas for central extended binomial
  coefficients, which 
are generalizations of binomial coefficients.
To do so, we relate the exact distribution of the sum of independent
discrete uniform random variables to the asymptotic distribution,
obtained from the Central Limit Theorem and a 
local limit variant.
\end{abstract}

\section{Stirling's formula and central binomial coefficients}
For a nonnegative integer $k$, Stirling's formula
\begin{align*}
  k! \sim \sqrt{2\pi k}\left(\frac{k}{e}\right)^k
\end{align*}
where $e$ is Euler's number, yields an approximation of the central
binomial coefficient $\binom{k}{k/2}$ using
$\binom{k}{m}=\frac{k!}{m!(k-m)!}$ as 
\begin{align*}
  \binom{k}{k/2} \sim \frac{2^{k+1}}{\sqrt{2\pi k}},
\end{align*}
where we write $a_k\sim b_k$ as short-hand for
$\lim_{k\rightarrow\infty}\frac{a_k}{b_k}=1$. In our current note, we
derive asymptotic formulas for central \emph{extended binomial}, or
\emph{polynomial}, coefficients (cf.\ \cite{Caiado:2007,Comtet:1970,Fahssi:2012}). These coefficients
appear in the \emph{extended binomial triangles} (which we also call
\emph{$(\ell+1)$-nomial, polynomial}, or \emph{multinomial triangles}
\cite{Fielder:1991}), which are generalizations
of binomial, or Pascal, triangles, 
in which entries in row $k$ are defined as coefficients of the
polynomial $(1+x+x^2+\cdots+x^\ell)^k$ for $\ell\ge 0$. Our derivation is
not based upon asymptotics of factorials, but upon the limiting
distribution of the sum of discrete uniform random
variables.\footnote{Throughout, we assume that all fractional values
  such as $x=\frac{k\ell}{2}$ are integral when used in the context of
  extended binomial coefficients. If this is not the case, then
  replace respective quantities with their \emph{floor}, $\lfloor
  x\rfloor$, the largest integer less than or equal to $x$.}

\section{Extended binomial triangles}
In generalization to binomial triangles, $(\ell+1)$-nomial triangles, for
$\ell\ge 0$, are defined in the following way. Starting with a $1$ in row
zero, construct an entry in row $k$, for $k\ge 1$, by adding the
overlying $(\ell+1)$ entries in row $(k-1)$ (some of these entries are
taken as zero if not defined); thereby, row $k$ has $(k\ell+1)$
entries. For example, the binomial ($\ell=1$), trinomial $(\ell=2)$, and
quadrinomial triangles $(\ell=3)$ start as follows, 

\begin{table}[!htb]
  \begin{center}
  \begin{tabular}{cccc}
    1 \\
    1 & 1 \\
    1 & 2 & 1\\
    1 & 3 & 3 & 1\\
  \end{tabular}
  \begin{tabular}{ccccccc}
    1 \\
    1 & 1 & 1\\
    1 & 2 & 3 & 2 & 1\\
    1 & 3 & 6 & 7 & 6 & 3 & 1\\
  \end{tabular}
  \begin{tabular}{cccccccccc}
    1 \\
    1 & 1 & 1 & 1\\
    1 & 2 & 3 & 4 & 3 & 2 & 1\\
    1 & 3 & 6 & 10 & 12 & 12 & 10 & 6 & 3 & 1\\
  \end{tabular}
  \end{center}
\end{table}
In the $(\ell+1)$-nomial triangle, the $n$th entry, for $0\le n\le k\ell$ in
row $k$, which we denote by $\binom{k}{n}_{\ell+1}$, has the following
interpretation. It is the coefficient of $x^n$ in the expansion of 
\begin{align}\label{eq:1}
  (1+x+x^2+\cdots+x^\ell)^k = \sum_{n=0}^{k\ell}\binom{k}{n}_{\ell+1}x^n.
\end{align}
It has been shown that $\binom{k}{n}_{\ell+1}$ denotes the number of
restricted \emph{integer compositions} (for a definition, see, e.g.,
\cite{Flajolet:2009} and many others) of the nonnegative integer $n$
with $k$ parts 
$\pi_1,\ldots,\pi_k$, each from the set $\set{0,1,\ldots,\ell}$
(cf.\ \cite{Eger:2013}), and allows the following representation,
\begin{align}\label{eq:2}
  \binom{k}{n}_{\ell+1} = \sum_{\substack{k_0\ge 0,\ldots,k_\ell\ge
      0\\ k_0+\cdots+k_\ell=k \\ 0\cdot k_0+1\cdot k_1+\cdots+\ell\cdot
      k_\ell=n}}\binom{k}{k_0,\ldots,k_\ell}, 
\end{align}
where $\binom{k}{k_0,\ldots,k_\ell}$ is a \emph{multinomial coefficient},
defined as $\frac{k!}{k_0!\cdots k_\ell!}$, for nonnegative integers
$k_0,\ldots,k_\ell$. We can verify representation \eqref{eq:2} by noting
that for real numbers $x_0,\ldots,x_\ell$, the multinomial theorem
(cf.\ \cite{Weisstein}) states that
\begin{align*}
  (x_0+x_1+\cdots+x_\ell)^k = \sum_{\substack{k_0\ge 0,\ldots,k_\ell\ge
      0\\ k_0+\cdots+k_\ell=k }}\binom{k}{k_0,\ldots,k_\ell}x_0^{k_0}\cdots
  x_\ell^{k_\ell}. 
\end{align*}
Thus, setting $x_i=x^i$ for $i=0,\ldots,\ell$,
\begin{align}\label{eq:3}
  (1+x^1+\cdots+x^\ell)^k = \sum_{\substack{k_0\ge 0,\ldots,k_\ell\ge
      0\\ k_0+\cdots+k_\ell=k }}\binom{k}{k_0,\ldots,k_\ell}x^{0\cdot
    k_0+\cdots+\ell\cdot k_\ell},
\end{align}
so that comparing coefficients of the right-hand sides of \eqref{eq:1}
and \eqref{eq:3} leads to \eqref{eq:2}.

\section{Generalized Stirling's approximation} 
Our strategy for deriving approximation formulas for central extended
binomial coefficients is as follows. First, we determine the
asymptotic distribution of the sum of discrete uniform variables,
which we easily find to be a normal distribution by the Central Limit
Theorem (CLT). Then, we determine the exact distribution, which turns
out to yield the normalized extended binomial coefficients
$\binom{k}{n}_{\ell+1}$. By relating the density of the asymptotic
distribution to the density of the exact distribution (e.g., via a
`local limit' argument), we obtain an extended binomial analgoue of
Stirling's approximation to central binomial coefficients.

\subsection{Step 1: Asymptotic distribution of the sum of discrete
  uniform variables}
Let $k$ be a positive integer and let $\ell$ be a nonnegative
integer. Let $X_j$, for $j=1,\ldots,k$, be identically and
independently distributed random draws from the discrete uniform
distribution on the set $\set{0,\ldots,\ell}$, and let $S_k$ be their
sum,
\begin{align*}
  S_k = \sum_{j=1}^k X_j.
\end{align*}
By standard moments of the uniform distribution, the mean and variance
of each $X_j$ are given by
\begin{align*}
  \mu = \Exp[X_j]=\frac{\ell}{2},\quad\text{and}\quad
  \sigma^2=\Var[X_j] = \frac{(\ell+1)^2-1}{12}.
\end{align*}
Hence, by independent and identical distribution of $X_1,\ldots,X_k$,
and application of the CLT, the random variable
$\sqrt{k}(\frac{S_k}{k}-\mu)$ converges, as $k\rightarrow\infty$, in
distribution to a normal $\mathcal{N}(0,\sigma^2)$ distributed random
variable. Recall that convergence in distribution precisely means that
the cumulative density function of $\sqrt{k}(\frac{S_k}{k}-\mu)$
converges pointwise to the cumulative density function of the
$\mathcal{N}(0,\sigma^2)$ distribution. 

\subsection{Step 2: Exact distribution of the sum of discrete uniform random variables}
We now determine exactly the probability that $S_k$ takes on the
integer value $n$, for $0\le n\le k\ell$. To do so, we consider
`isomorphic copies' $\tilde{X}_j$ of $X_j$, which are independently
and identically \emph{multinomially} distributed with probabilities
$p_0=\cdots = p_{\ell}=\frac{1}{\ell+1}$ of types $0$ to $\ell$. Each
$\tilde{X}_j=(A_0,\ldots,A_\ell)$ is vector-valued, with
$P[\tilde{X}_j=(a_0,\ldots,a_\ell)]=\frac{1}{\ell+1}$ for nonnegative
integers $a_s$, with $a_0+\cdots+a_\ell=1$, where $A_s$ denotes the
number of times an event of type $s$, for $s=0,\ldots,\ell$,
occurs. Then, the sum $\tilde{S}_k=\tilde{X}_1+\cdots+\tilde{X}_k$ has
the interpretation of representing the event of drawing with
replacement $k$ balls of $(\ell+1)$ different types from a bag, where
the probability of drawing type $s=0,\ldots,\ell$ is
$\frac{1}{\ell+1}$. Thus, by the standard interpretation of the
multinomial distribution, $\tilde{S}_k$ has density
\begin{align*}
  P[\tilde{S}_k=(a_0,\ldots,a_\ell)] = P[A_0=a_0,\ldots,A_\ell=a_\ell]
  = \binom{k}{a_0,\ldots,a_\ell}\left(\frac{1}{\ell+1}\right)^k,
\end{align*}
where $a_0+\cdots+a_\ell=k$ for nonnegative integers
$a_0,\ldots,a_\ell$. Then, if $\tilde{S}_k=(a_0,\ldots,a_\ell)$,
$S_k$, the variable corresponding to $\tilde{S}_k$, represents the
integer $0\cdot a_0+\cdots+\ell\cdot a_\ell$. Thus, for $n$ such that
$0\le n\le k\ell$,
\begin{align*}
  P[S_k=n]=\sum_{\substack{a_0\ge 0,\ldots,a_\ell\ge
      0\\ a_0+\cdots+a_\ell=k \\ 0\cdot a_0+\cdots+\ell\cdot
      a_\ell=n}} P[\tilde{S}_k=(a_0,\ldots,a_\ell)]
  =\left(\frac{1}{\ell+1}\right)^k\binom{k}{n}_{\ell+1},
\end{align*}
using representation \eqref{eq:2}.

An arguably more straightfoward derivation of the exact distribution
of $S_k$, making use of probability generating functions (pgfs), can
be given by noting that the pgf $G_{X_j}(x)=\sum_{n\ge 0}P[X_j=n]x^n$
of each $X_j$ is given by
\begin{align*}
  G_{X_j}(x)=\frac{1}{\ell+1}\sum_{n=0}^\ell x^n.
\end{align*}
By independence of $X_1,\ldots,X_k$, the pgf of $S_k$ is hence given as,
\begin{align*}
  G_{S_k}(x) = G_{X_1}(x)\cdots
  G_{X_k}(x)=\left(\frac{1}{\ell+1}\right)^k\left(\sum_{n=0}^\ell
  x^n\right)^k =
  \left(\frac{1}{\ell+1}\right)^k\sum_{n=0}^{k\ell}\binom{k}{n}_{\ell+1}x^n.
\end{align*}
Thus,
\begin{align*}
  P[S_k=n]=\frac{G_{S_k}^{(n)}(0)}{n!} =
  \left(\frac{1}{\ell+1}\right)^k\frac{n!}{n!}\binom{k}{n}_{\ell+1}=\left(\frac{1}{\ell+1}\right)^k\binom{k}{n}_{\ell+1},
\end{align*}
where we denote by $G_X^{(n)}(0)$ the $n$th derivative of $G_X$,
evaluated at zero. 

\subsection{Step 3: Local limit theorem}
\begin{figure}[!htb]
  \begin{center}
  \input{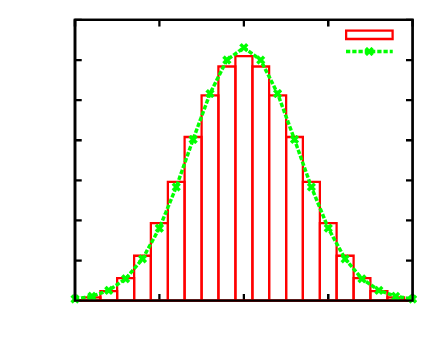}
  \input{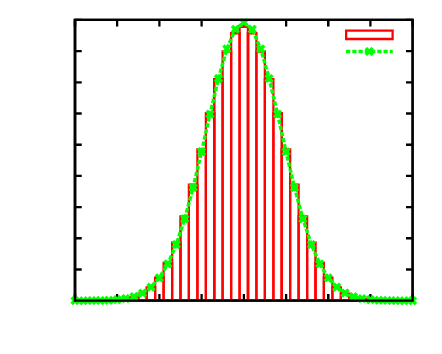}
  \input{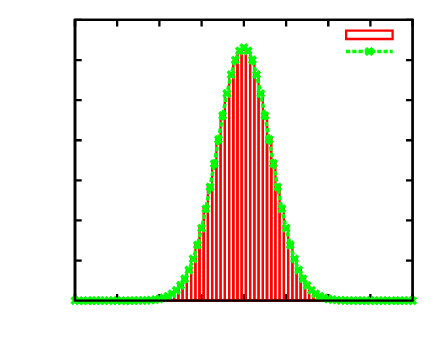}
  \caption{Distributions $P[S_k=n]$ for $k=5,10,20$, for $l=4$ fixed, and normal approximations.}
  \label{fig:1}
  \end{center}
\end{figure}

To derive an asymptotic formula for $\binom{k}{n}_{\ell+1}$, we
would like to make use of the results derived in Steps 1 and 2
above. Ideally, we would 
like to equate the probability density function of the asymptotic
normal dstribution of $S_k$
with the exact distribution. However, as mentioned, convergence in distribution,
as assured by the CLT, only guarantees pointwise convergence of cumulative density
functions. On the contrary, `local limit theorems' describe how the
probability density function of a sum of random variables approaches
the normal density function. 
For integer-valued random variables (also called lattice or
arithmetical distributions), 
Gnedenko and Kolmogorov \cite{Gnedenko:1968} provide the following result.
\begin{theorem}\label{th:1}
If $X_1,X_2,\ldots$ are independent lattice random variables with
identical distribution with finite mean $\mu$ and variance $\sigma^2$,
such that the greatest common divisor of the differences of all the
values of $X_j$ taken with positive probability is $1$, then
\begin{align*}
  \left|\sqrt{k}\sigma
  P[S_k=n]-\frac{1}{\sqrt{2\pi}}e^{-\frac{(n-k\mu)^2}{2\sigma^2k}}\right|\rightarrow 0
\end{align*}
uniformly in $n$ as $k\rightarrow\infty$.
\end{theorem}
Since in our situation, the set of values of each $X_j$ taken with
positive probability is $\set{0,\ldots,\ell}$, the greatest common
divisor of the differences is clearly $1$. Thus, all assumptions of
Theorem \ref{th:1} are satisfied in our case, and, hence, also
its consequences hold. Therefore, the following approximation is
suggested for large $k$:
\begin{align}\label{eq:approx}
  \sqrt{k}\sigma P[S_k=n]\sim \frac{1}{\sqrt{2\pi}}e^{-\frac{(n-k\mu)^2}{2\sigma^2k}}.
\end{align}
For $n=k\mu=k\ell/2$, the argument to the exponential function is
zero, and thus
\begin{align*}
  \sqrt{k}\sigma P[S_k=k\ell/2]\sim
  \frac{1}{\sqrt{2\pi}},\quad\text{or equivalently},\quad
  P[S_k=k\ell/2]\sim \frac{1}{\sqrt{2\pi\sigma^2k}}.
\end{align*}
Using the exact form for $P[S_k=n]$ from Step 2 above, we hence have,
bringing the normalizing term $(\ell+1)^k$ to the right-hand side,
\begin{align}\label{eq:main}
  \binom{k}{\frac{k\ell}{2}}_{\ell+1}\sim \frac{(\ell+1)^k}{\sqrt{2\pi
  k\frac{(\ell+1)^2-1}{12}}}.
\end{align}
For example, for $\ell=1$, Pascal's case, $\ell=2$, $\ell=3$, and
$\ell=4$, we therefore have the approximations
\begin{align*}
  \binom{k}{\frac{k}{2}}\sim \frac{2^{k+1}}{\sqrt{2\pi k}},\quad 
  \binom{k}{k}_3 \sim \frac{3^k}{\sqrt{\frac{4}{3}\pi k}},\quad
  \binom{k}{\frac{3}{2}k}_4 \sim \frac{4^k}{\sqrt{\frac{5}{2}\pi
      k}},\quad\text{and}\quad
  \binom{k}{2k}_5\sim \frac{5^k}{2\sqrt{\pi k}}.
\end{align*}
In Figure \ref{fig:1}, we show for $\ell=4$ the distributions
$P[S_k=n]$ for $k=5,10,20$, and their respective normal
approximations. There, we can see the local limit theorem `at work':
The exact density function apparently approaches, pointwise, the
normal density function. 

\section{Discussion}
Although extended binomial coefficients, together with their connection to the sum of discrete uniform random variables, go back at least to De
Moivre's Doctrine of Chances \cite{DeMoivre:1967} and to Euler's
\cite{Euler:1801} analytical study of the coefficients 
of polynomial \eqref{eq:1}, the mathematics community has apparently
more or less ignored  
their systematic study, except for a few recent publications such as
\cite{Bollinger:1990,Caiado:2007,Eger:2013,Fahssi:2012,Fielder:1991}. 
Next, using the CLT (or a local limit variant) to deduce asymptotics
of mathematical objects has been suggested, for example, by Walsh
\cite{Walsh:1995}, who derives Stirling's 
formula for factorials by equating the distribution of the sum of
Poisson distributed 
random variables with the normal density. Finally, the asymptotics of
both the central 
binomial ($\ell=1$) as well as the central trinomial coefficients
($\ell=2$) seem to be known 
(e.g.~\cite{Fahssi:2012,OEIS}), while the general formula \eqref{eq:main} is, to the best of our knowledge, novel.
However, Ratsaby \cite{Ratsaby:2008} derives our general result \eqref{eq:approx}, as an estimate of the number
of restricted integer compositions, by application of Cauchy's
coefficient formula to 
the polynomial \eqref{eq:1} and computation of the resulting integral
by Laplace's method for 
evaluation of integrals. A historical perspective of local versus
central limit theorem 
is provided by McDonald \cite{MacDonald:2005}.

\bibliographystyle{plain}
\bibliography{lit}

\end{document}